\def\reporttype{1}

\def\INRIAreport{1}
\ifnum\reporttype<0\stop\fi
\ifnum\reporttype>1\stop\fi

\documentclass[11pt]{article}
\usepackage{makeidx}
\usepackage{amssymb}
\usepackage{amsmath}
\usepackage[dvips]{graphicx}
\usepackage{color}
\usepackage{url}
\usepackage{multicol}

\ifnum\reporttype=\INRIAreport
\usepackage{RR}
\fi

\ifnum\reporttype=\INRIAreport
\else
\fi

\def\withRCSversion{0}
\def\withaddnotes{0}
\def\withenglishnotes{0}
\def\withffoonotes{0}

\newcommand{\addnote}[1]{\ifnum\withaddnotes=0\else{\footnote{#1}}\fi}
\newcommand{\englishnote}[1]{\ifnum\withenglishnotes=0%
 \else{\footnote{{\em English\/}: #1}}\fi}
\newcommand{\ffootnote}[1]%
 {\ifnum\withffoonotes=0\else{\footnote{#1}}\fi}

\newlength{\minipagewidth}

\newenvironment{accolade}{
 \left\{\begin{array}
 }{\end{array}\right.}
 {\unskip\nobreak\hfill\penalty80\hskip1em\hbox{}\nobreak\hfill$\Box$}

\def\defrcsversion$#1: modulopt.tex,v #2 #3/#4#5/#6#7 #8 #9 $%
 {\def\doubleday{#6}\def\doublemonth{#4}%
  \def\rcsversion{#2 -- %
  \ifnum\doubleday=0\else{#6}\fi#7/%
  \ifnum\doublemonth=0\else{#4}\fi#5/#3 -- #8}}
\defrcsversion$Id: modulopt.tex,v 1.8 2007/02/21 14:20:55 gilbert Exp $

\def\mat#1#2{%
  \if#1<%
    \if#2=\preccurlyeq\else\prec#2\fi%
  \else\if#1>%
         \if#2=\succcurlyeq\else\succ#2\fi%
       \else#1#2%
       \fi%
  \fi}

\def\<#1,#2>{\langle #1,#2\rangle}

\newcommand{\AND}{\quad\mbox{and}\quad}
\newcommand{\PROB}{\textit{\texttt{prob}}}
\newcommand{\RR}{{\mathbb R}}
\newcommand{\RRb}{\gbar{\RR}}
\newcommand{\RRm}{\RR^m}
\newcommand{\RRmE}{\RR^{m_E}}
\newcommand{\RRmI}{\RR^{m_I}}
\newcommand{\RRn}{\RR^n}
\newcommand{\SOLV}{\textit{\texttt{solv}}}
\newcommand{\T}{^{\mskip-1mu\top\mskip-2mu}}

\newcommand{\abr}{\allowbreak}
\newcommand{\gbar}[1]{\newbox\charbox\setbox\charbox=\hbox{$#1$}%
 \vbox{\vbox{\hrule width0.05\wd\charbox height0pt%
             \hrule width0.9\wd\charbox height0.3pt}%
       \nointerlineskip\kern0.12em\box\charbox}}
\newcommand{\lb}{\bar{\lm}}
\newcommand{\libopt}{\texttt{libopt}}
\newcommand{\lm}{\lambda} 
\newcommand{\sgn}{\operatorname{sgn}} 
\newcommand{\union}{\cup} 
\newcommand{\xb}{\bar{x}}

\makeindex

\index{nonsmooth|see{function/nonsmooth}}
\index{nondifferentiable|see{function/nonsmooth}}
\index{modulopttoys@\texttt{modulopttoys}|see{collection}}

\newbox{\mybox}

\makeatletter
  {\end{multicols}}
\makeatother

\newcommand{\theabstract}{%
This note describes how the optimization problems of the Modulopt
collection are organized within the Libopt environment. It is aimed at
being a guide for using and enriching this collection in this
environment.
}

\ifnum\reporttype=\INRIAreport\relax
\RRtitle{%
Organisation de la collection de probl\`emes d'optimisation Modulopt
dans l'environnement Libopt\\-- Version 1.0 --
}
\RRetitle{%
Organization of the Modulopt collection of optimization problems in the
Libopt environment\\-- Version 1.0 --
}
\titlehead{Organization of Modulopt in Libopt}
\RRauthor{%
  J.\ Charles {\sc Gilbert}\/%
  \thanks{INRIA Rocquencourt, projet Estime, BP~105, 78153~Le Chesnay
  Cedex, France\,;
  e-mail\,: {\tt Jean}\abr{\tt -}\abr{\tt Charles.}\abr%
            {\tt Gilbert@}\abr{\tt inria.fr}.}
}
\authorhead{J.\ Ch.\ Gilbert}
\RRdate{21 f\'evrier 2007}
\RRresume{%
Cette note d\'ecrit comment les probl\`emes d'optimisation de la
collection Modulopt sont organis\'es dans l'environnement Libopt. Elle
a pour but de servir de guide pour utiliser et enrichir cette
collection dans cet environnement.
}
\RRabstract{\theabstract}
\RRmotcle{%
collection de probl\`emes --
environnement de test --
\'evaluation de performance --
Libopt --
Modulopt --
optimisation
}
\RRkeyword{%
benchmarking --
collection of problems --
Libopt --
Modulopt --
optimization --
testing environment
}
\RRprojet{Estime}
\RRtheme{\THNum}
\URRocq
\fi

\begin{document}

\ifnum\reporttype=\INRIAreport\else
\begin{center}
{\LARGE\bf
Organization of the Modulopt collection\\
of problems in the Libopt environment\\}
\vspace{2ex}\rm\normalsize
Version 1.0~~~(\today)
\\
\vspace{2ex}\rm\normalsize
Jean Charles {\sc Gilbert}\/\footnote{INRIA-Rocquencourt, BP~105,
F-78153~Le Chesnay Cedex (France); e-mail: \texttt{Jean-Charles.}\abr
\texttt{Gilbert@}\abr \texttt{inria.fr}.}

\ifnum\withRCSversion=0\else{%
\vspace*{.2cm}\\
{\bf Version \rcsversion}}\fi
\end{center}

\renewcommand{\thefootnote}{\fnsymbol{footnote}}
\setcounter{footnote}{1}

\renewcommand{\thefootnote}{\arabic{footnote}}
\setcounter{footnote}{1}

\begin{abstract}
\noindent
\theabstract
\end{abstract}
\fi


\ifnum\reporttype=\INRIAreport\relax
\makeRT
\fi


\section{The problems of the Modulopt collection}

In the Libopt terminology~\cite{gilbert-jonsson-2007}, a {\em
collection\/}\index{collection} refers to a set of problems sharing
some common features, such as their mathematical structure, coding
language, audience, etc. In this note, we describe the installation of
the Modulopt collection \cite{lemarechal-1980} in the Libopt
environment. The features of the Modulopt problems, from the Libopt
viewpoint, are the following:
\begin{list}{{\small$\bullet$}}{\topsep=1.5ex\parsep=0.0ex\itemsep=0.0ex
 \settowidth{\labelwidth}{{\small$\bullet$}}
 \labelsep=0.5em
 \leftmargin=\labelwidth
 \addtolength{\leftmargin}{\labelsep}
 \addtolength{\leftmargin}{\parindent}
}
\item
they have an optimization nature and can be written in the
form~\eqref{p} below;
\item
they can be smooth or nonsmooth;
\item
they are written in Fortran 90/95;
\item
they are issued from various application areas in scientific or
industrial computing;
\item
they can be freely distributed.
\end{list}
The collection has a companion one, named
``\texttt{modulopttoys}''\index{collection!modulopttoys@\texttt{modulopttoys}},
which has the same features, except that the problems have an academic
nature. In Libopt, these collections rub shoulders with the {\em CUTEr
collection\/}\index{collection!CUTEr}~\cite{bongartz-conn-gould-toint-1995,
gould-orban-toint-2003}.

The Modulopt collection contains nonlinear optimization problems coming
from various application areas. The optimization problems are supposed
to be written in the following form
\begin{equation}
\label{p}
(P)\quad
\begin{accolade}{l}
\min\;f(x)\\
l_B\leq x\leq u_B\\
l_I\leq c_I(x)\leq u_I\\
c_E(x)=0,
\end{accolade}
\end{equation}
where $f:\RRn\to\RR$, $c_I:\RRn\to\RRmI$, $c_E:\RRn\to\RRmE$, $l_B$,
$u_B\in\RRb^n$, and $l_I$, $u_I\in\RRb^{m_I}$
($\RRb=\RR\union\{-\infty,+\infty\}$). Actually $B$ is the set of
indices $\{1,\ldots,n\}$ and $I$ is another set of indices with $m_I$
elements. We write $l:=(l_B,l_I)\in\RRb^n\times\RRb^{m_I}$ and
$u:=(u_B,u_I)\in\RRb^n\times\RRb^{m_I}$. It is assumed that $l<u$,
meaning that $l_i<u_i$, for all $i\in B\union I$. For making the
notation compact, we note
$$
c_B(x):=x,
\quad
c(x):=(c_B(x),c_I(x),c_E(x)),
\AND
m:=n+m_I+m_E.
$$
The Jacobian matrices of $c_I$ and $c_E$ at $x\in\RRn$ are also denoted
by
$$
A_I(x):=c'_I(x)
\AND
A_E(x):=c'_E(x).
$$
We also introduce the nondifferentiable operator
$(\cdot)^\#:\RRm\to\RRm$ defined by
$$
v^\#=
\begin{pmatrix}
\max(0,l_B-v_B,v_B-u_B)\\
\max(0,l_I-v_I,v_I-u_I)\\
v_E
\end{pmatrix},
$$
so that $x$ is feasible for $(P)$ if and only if $c(x)^\#=0$.

The {\em Lagrangian\/}\index{Lagrangian} of problem $(P)$ is the
function $\ell:\RRn\times\RRm\to\RR$ defined at $(x,\lm)$~by
\begin{equation}
\label{lagrangian}
\ell(x,\lm) = f(x) + \lm\T c(x).
\end{equation}
Note that we take a single multiplier for two constraints present in the
bound constraints $l_i\leq c_i(x)\leq u_i$, knowing that $l_i<u_i$
implies that at least one of the multipliers associated with $l_i\leq
c_i(x)$ and $c_i(x)\leq u_i$ is zero. The {\em optimality
conditions\/}\index{optimality conditions} at $\xb$ read for some
optimal multiplier~$\lb$\ffootnote{See \cite{dolan-more-munson-2006}
for a more convenient stopping test.}:
\begin{equation}
\label{kkt}
\begin{accolade}{l}
\nabla f(\xb)+c'(\xb)\T\lb=0\\
c(\xb)^\#=0\\
i\in B\union I,~l_i<c_i(\xb)<u_i
\quad\Longrightarrow\quad\lb_i=0\\
i\in B\union I,~l_i=c_i(\xb)\hphantom{{}<u_i}
\quad\Longrightarrow\quad\lb_i\leq0\\
i\in B\union I,~\hphantom{l_i<{}}c_i(\xb)=u_i
\quad\Longrightarrow\quad\lb_i\geq0.
\end{accolade}
\end{equation}

\section{Running a Modulopt problem}

\subsection{Notation and relevant directories}

We use the following typographic conventions. The \texttt{typewriter
font} is used for a text that has to be typed literally and for the
name of files and directories that exist as such (without making
substitutions). In the same circumstances, a generic word, which has to
be substituted by an actual word depending on the context, is written in
\texttt{\textit{italic typewriter font}}.

Here are some directories of the Libopt hierarchy that will intervene
continually in this note. Other important directories and files
introduced in this note are listed in
section~\ref{s:directories-files}. The main directories~are
\begin{list}{{\small$\bullet$}}{\topsep=1.0ex\parsep=0.5ex\itemsep=0.5ex
 \settowidth{\labelwidth}{{\small$\bullet$}}
 \labelsep=0.5em
 \leftmargin=\labelwidth
 \addtolength{\leftmargin}{\labelsep}
 \addtolength{\leftmargin}{\parindent}
 \rightmargin=\parindent
}
\item
\texttt{\$LIBOPT\_DIR}\\%
\index{environment variable!\texttt{\$LIBOPT\_DIR}}%
\index{libopt!environment variable}%
is the environment variable that specifies the {\em root
directory\/}\index{directory!libopt root@\texttt{libopt}
root}\index{libopt!root directory} of the Libopt hierarchy,
\item
\texttt{\$LIBOPT\_DIR/collections/modulopt}\\
is the {\em root directory of the Modulopt
collection\/}\index{directory!Modulopt root} in the Libopt environment,
\item
\texttt{\$LIBOPT\_DIR/collections/modulopt/probs}\\
is the directory that has a sub-directory for each of the problems of
the Modulopt collection installed in the Libopt environment,
\item
\texttt{\$LIBOPT\_DIR/solvers}\\
is the {\em root directory of the solvers\/}\index{directory!solver
root} installed in the Libopt environment.
\end{list}

\subsection{The \texttt{runopt} script}

\index{runopt@\texttt{runopt} (script)|(}

The simplest way of running a single Modulopt problem in the
Libopt environment is by typing (`\texttt{\%}' is the Unix/Linux prompt)
\begin{quote}
\texttt{%
\% echo "\SOLV\ modulopt \PROB" | runopt},
\end{quote}
where, here and below,
\begin{list}{{\small$\bullet$}}{\topsep=1.0ex\parsep=0.5ex\itemsep=0.5ex
 \settowidth{\labelwidth}{{\small$\bullet$}}
 \labelsep=0.5em
 \leftmargin=\labelwidth
 \addtolength{\leftmargin}{\labelsep}
 \addtolength{\leftmargin}{\parindent}
 \rightmargin=\parindent
}
\item
\SOLV\index{solv@\SOLV} stands for the name of a solver installed
in the Libopt environment, one of those listed in
\begin{quote}
\texttt{%
\$LIBOPT\_DIR/solvers/solvers.lst};
\end{quote}
actually, the solver \SOLV\ must have been prepared to run Modulopt
problems, otherwise this command will not be understood by the Libopt
environment; this subject is discussed in
section~\ref{s:preparing-solver};
\item
\PROB\index{prob@\PROB} stands for the name of a Modulopt
problem currently available in the Libopt environment, one in the list
\begin{quote}
\texttt{%
\$LIBOPT\_DIR/collections/modulopt/all.lst}.
\end{quote}
\end{list}
By this command the optimization code \SOLV\ is used to solve the
Modulopt optimization problem \PROB. Of course \SOLV\ has to be able to
solve a problem with the features of \PROB\ (for example, a solver for
unconstrained optimization problems is unable to solve problems with
constraints). The code \SOLV\ keeps in the file
\begin{quote}
\texttt{\$LIBOPT\_DIR/solvers/\SOLV/modulopt/all.lst}
\end{quote}
the list of the Modulopt problems that it can structurally solve.

See \cite{gilbert-jonsson-2007} or the manual page of \texttt{runopt}
to learn how to run a group of problems with a given solver, using a
single command line or a file describing what has to be done.

\index{runopt@\texttt{runopt} (script)|)}

The directory where the \texttt{runopt} script given above is typed is
called the {\em working directory}\index{directory!working}. The Libopt
scripts take care that this directory is not in the Libopt hierarchy.
If this were the case, there would be a danger of incurable
destruction. Indeed, a script like
\texttt{runopt}\index{command!runopt@\texttt{runopt}} generally removes
several files from the working directory after a problem has been
solved.

\subsection{The \texttt{\SOLV\_modulopt} script}

\index{solver\_modulopt@\texttt{\SOLV\_modulopt} (script)|(}

By decoding its standard input ``\texttt{\SOLV\ modulopt \PROB}'',
the \texttt{runopt} script above knows that it has to launch the
following command:
\begin{quote}
\texttt{%
\$LIBOPT\_DIR/solvers/\SOLV/modulopt/\SOLV\_modulopt \PROB}.
\end{quote}
In the standard distribution, \texttt{\SOLV\_modulopt} is a Perl
script, but nothing imposes that such a language be used. Such a script
has to be written for each solver (actually for each solver-collection
pair). Section~\ref{s:preparing-solver} explains how to do this. For
the while, it is enough to know that it is decomposed in the following
main steps.
\begin{list}{{\small$\bullet$}}{\topsep=1.0ex\parsep=0.5ex\itemsep=0.5ex
 \settowidth{\labelwidth}{{\small$\bullet$}}
 \labelsep=0.5em
 \leftmargin=\labelwidth
 \addtolength{\leftmargin}{\labelsep}
 \addtolength{\leftmargin}{\parindent}
 \rightmargin=\parindent
}
%
\item
The environment variables
\begin{quote}
\texttt{\$MODULOPT\_PROB}\index{environment
variable!\texttt{\$MODULOPT\_PROB}} and
\texttt{\$WORKING\_DIR}\index{environment
variable!\texttt{\$WORKING\_DIR}}
\end{quote}
are respectively set to \PROB\ and to the {\em working
directory}\index{directory!working}, so that these variables can be
used in the makefiles mentioned below. Actually,
\texttt{\$WORKING\_}\abr \texttt{DIR} is probably useless since all the
Unix/Linux commands in the scripts or makefiles are executed from the
working directory (there is no change of directory made in them).

\item
Then the target \PROB\ of the following makefile is executed
\begin{quote}
\texttt{%
\$LIBOPT\_DIR/collections/modulopt/probs/\PROB/Makefile}.
\end{quote}
The aim of this target is to make symbolic links in the working
directory to the source and data files in the \PROB\ directory
\begin{quote}
\texttt{%
\$LIBOPT\_DIR/collections/modulopt/probs/\PROB},
\end{quote}
to produce an archive named \texttt{\PROB.a} in the working directory,
which contains the problem object files allowing the execution of the
problem, and finally to remove the now useless symbolic links from the
working directory.

\item
Next, the Perl script executes the target
\texttt{\SOLV\_modulopt\_main} of the following makefile
\begin{quote}
\texttt{%
\$LIBOPT\_DIR/solvers/\SOLV/modulopt/Makefile}.
\end{quote}
Its aim is to make a symbolic link in the working directory to the
source file
\begin{quote}
\texttt{%
\$LIBOPT\_DIR/solvers/\SOLV/modulopt/\SOLV\_modulopt\_main.f90}
\end{quote}
of the main program, to compile it and link it with the archive
\texttt{\PROB.a} of the Modulopt problem previously generated. This
produces the executable file \texttt{\SOLV\_modulopt\_main} in the
working directory. Then the target removes from the working directory
the now useless symbolic link \texttt{\SOLV\_}\abr
\texttt{modulopt\_}\abr \texttt{main.f90} and file \texttt{\SOLV\_}\abr
\texttt{modulopt\_}\abr \texttt{main.o}.

\item
The program
\begin{quote}
\texttt{%
\SOLV\_modulopt\_main}
\end{quote}
is then executed in the working directory. This one solves the problem
\PROB\ with the solver \SOLV.

\item
Some cleaning is then done in the working directory:
\texttt{\SOLV\_}\abr \texttt{modulopt\_}\abr \texttt{main} is removed
(probably with other files, depending on the solver) and the target
\texttt{\PROB\_}\abr \texttt{clean} of the following makefile is
executed:
\begin{quote}
\texttt{%
\$LIBOPT\_DIR/collections/modulopt/probs/\PROB/Makefile}
\end{quote}
Its aim is to remove from the working directory, the files related to
the problem just solved, typically the archive \texttt{\PROB.a} and
the \PROB\ data files.
\end{list}

\index{solver\_modulopt@\texttt{\SOLV\_modulopt} (script)|)}

\section{Introducing a new problem in the Modulopt collection}

\subsection{Overview}
\label{s:new-problem}

Suppose that we want to introduce a new problem in the Modulopt
collection and that this one is called
\begin{quote}
\PROB\index{prob@\PROB}.
\end{quote}
The description of the \texttt{runopt} script above shows that, one has
to proceed as follows.
\begin{list}{{\small$\bullet$}}{\topsep=1.0ex\parsep=0.5ex\itemsep=0.5ex
 \settowidth{\labelwidth}{{\small$\bullet$}}
 \labelsep=0.5em
 \leftmargin=\labelwidth
 \addtolength{\leftmargin}{\labelsep}
 \addtolength{\leftmargin}{\parindent}
 \rightmargin=\parindent
}
\item
Insert the name \PROB\ in the list of Modulopt problems
\begin{quote}
\begin{verbatim}
$LIBOPT_DIR/collections/modulopt/all.lst
\end{verbatim}
\end{quote}
(and possibly in other lists in the same directory, such as the one
related to unconstrained problems \texttt{unc.lst}, quadratic problems
\texttt{quad.lst}, etc, as well as the list of typical problems of the
Modulopt collection \texttt{default.lst}, if this is appropriate). This
is an \texttt{ascii} file. An alpha-numeric order has been adopted, but
this feature is not taken into account by the Libopt scripts. {\em
Comments\/}\index{list!comment in a --} are possible; they start from
the character `\texttt{\#}' up to the end of the~line.

\item
If a solver called \SOLV\ is able to solve a problem like \PROB, it may
be appropriate to insert the name \PROB\ in one or more files among
\begin{quote}
\texttt{%
\$LIBOPT\_DIR/solvers/\SOLV/modulopt/*.lst}.
\end{quote}
This assumes that the directory
\texttt{\$LIBOPT\_DIR/solvers/\SOLV/modulopt} exists and that the
solver has been prepared to solve problems from the Modulopt collection
(see section~\ref{s:preparing-solver} to know how to do this).

\item
Create the directory
\begin{quote}
\texttt{%
\$LIBOPT\_DIR/collections/modulopt/probs/\textit{prob}},
\end{quote}
and put in that directory, all the files that define the problem
\PROB: source files, header files (if appropriate),
and data files (if appropriate). This is further described in
section~\ref{s:test-problem} below.

\item
Create the makefile
\begin{quote}
\texttt{%
\$LIBOPT\_DIR/collections/modulopt/probs/\textit{prob}/Makefile},
\end{quote}
with the following two targets:
\begin{list}{{\small$\circ$}}{\topsep=1.0ex\parsep=0.5ex\itemsep=0.5ex
 \settowidth{\labelwidth}{{\small$\circ$}}
 \labelsep=0.5em
 \leftmargin=\labelwidth
 \addtolength{\leftmargin}{\labelsep}
 \addtolength{\leftmargin}{\parindent}
 \rightmargin=\parindent
}
\item
\PROB, which specifies how to obtain in the working directory an
archive with all the object files defining \PROB\ and which makes
symbolic links in the working directory to the data files required to
solve the problem;
\item
\texttt{\textit{prob}\_clean}, which specifies which files has to be
removed from the working directory after having solved \PROB.
\end{list}
This is further discussed in section~\ref{s:makefile-test-problem}
below.
\end{list}

\subsection{The subroutines defining a Modulopt problem}
\label{s:test-problem}

In principle, the problem can be described in any compiled language,
provided the binary files can be gathered into an archive. Below, we
assume that the problem is written in \texttt{Fortran 95}.

The {\em problem-independent\/} makefile
\begin{quote}
\texttt{%
\$LIBOPT\_DIR/solvers/\SOLV/modulopt/Makefile}
\end{quote}
assumes that the problem to execute is in the archive
\texttt{\textit{prob}.a} in the working directory. On the other hand,
the {\em problem-independent\/} main program \texttt{\SOLV\_}\abr
\texttt{modulopt\_}\abr \texttt{main} assumes that the archive
\texttt{\textit{prob}.a} contains seven subroutines:
\texttt{dim}\-\texttt{opt}, \texttt{init}\-\texttt{opt},
\texttt{simul}\-\texttt{opt}, \texttt{postopt}, \texttt{inprodopt},
\texttt{ctonbopt}, and \texttt{ctcabopt}, which are described below.

In the description of the subroutine arguments, an argument tagged with
(I) means that it is an {\em input\/} variable, which has to be
initialized before calling the subroutine; an argument tagged with (O)
means that it is an {\em output\/} variable, which only has a meaning
on return from the subroutine; and an argument tagged with (IO) is an
{\em input-output\/} argument, which has to be initialized and which
has a meaning after the call to the subroutine. Arguments of the type
(O) and (IO) are generally modified by the subroutine and therefore
{\em should not be Fortran constants\/}!

\subsubsection*{The subroutine \texttt{dimopt}}

\index{subroutine!dimopt@\texttt{dimopt}|(}

The subroutine \texttt{dimopt} is called by the main program
\texttt{\SOLV\_}\abr \texttt{modulopt\_}\abr \texttt{main} to get the
dimensions of the problem. In \texttt{Fortran 95}, it has the following
calling structure:
\begin{quote}
\fbox{$
\setlength{\minipagewidth}{\linewidth}
\addtolength{\minipagewidth}{-6.8pt}
\begin{minipage}[t]{\minipagewidth}
\begin{list}{}{\topsep=0pt\partopsep=0pt\parsep=0.5\parskip\itemsep=0.5\parskip
 \settowidth{\labelwidth}{\texttt{subroutine dimopt (}}
 \labelsep=0ex
 \leftmargin=\labelwidth \addtolength{\leftmargin}{\parindent}
 \rightmargin=\parindent
 \raggedright
}
\item[\texttt{subroutine dimopt (}]
\texttt{n,} \texttt{mi,} \texttt{me,} \texttt{nizs,} \texttt{nrzs,}
\texttt{ndzs)}
\end{list}
\end{minipage}
$}
\end{quote}
\begin{list}{}{\topsep=0pt\partopsep=0pt\parsep=0.5\parskip\itemsep=0.5\parskip
 \listparindent=\parindent
 \labelwidth=\parindent
 \labelsep=1ex
 \addtolength{\leftmargin}{-1ex}
}
\item[\texttt{n} (O):]
positive \texttt{integer} variable. This is the number $n$ of variables
to optimize in the problem, those denoted $x=(x_1,\ldots,x_n)$ in
\eqref{p}.
\item[\texttt{mi} (O):]
nonnegative \texttt{integer} variable. This is the number $m_I$ of
nonlinear inequality constraints, of the form $l_i\leq c_i(x)\leq u_i$
($i=1,\ldots,m_I$), for some nonlinear functions
$c_i:\RRn\to\nobreak\RR$.
\item[\texttt{me} (O):]
nonnegative \texttt{integer} variable. This is the number $m_E$ of
nonlinear equality constraints, of the form $c_i(x)=0$
($i=1,\ldots,m_E$), for some nonlinear functions $c_i:\RRn\to\RR$.
\item[\texttt{nizs} (O), \texttt{nrzs} (O), \texttt{ndzs} (O):]
positive \texttt{integer} variables. These are the dimensions of the
variables \texttt{izs}, \texttt{rzs}, and \texttt{dzs} (respectively),
which are \texttt{integer}, \texttt{real}, and \texttt{double
precision} working zones for the Modulopt problem. The solvers must not
affect their content. The main program \texttt{\SOLV\_}\abr
\texttt{modulopt\_}\abr \texttt{main} associated with the solver \SOLV\
must allocate memory for the variables \texttt{izs}, \texttt{rzs}, and
\texttt{dzs} just after having called \texttt{dimopt}, see
section~\ref{s:preparing-solver}, point~4.1 on
page~\pageref{point-4.1}. This implies that using \texttt{Fortran}
\texttt{77} is not an appropriate language for writing the main program
\texttt{\SOLV\_}\abr \texttt{modulopt\_}\abr \texttt{main}. Note that
the value of \texttt{nizs}, \texttt{nrzs}, and \texttt{ndzs} can be
zero. \end{list}

\index{subroutine!dimopt@\texttt{dimopt}|)}

\subsubsection*{The subroutine \texttt{initopt}}

\index{subroutine!initopt@\texttt{initopt}|(}

The subroutine \texttt{initopt} is called to initialize the problem. In
\texttt{Fortran 95}, it has the following calling structure:
\begin{quote}
\fbox{$
\setlength{\minipagewidth}{\linewidth}
\addtolength{\minipagewidth}{-6.8pt}
\begin{minipage}[t]{\minipagewidth}
\begin{list}{}{\topsep=0pt\partopsep=0pt\parsep=0.5\parskip\itemsep=0.5\parskip
 \settowidth{\labelwidth}{\texttt{subroutine initopt (}}
 \labelsep=0ex
 \leftmargin=\labelwidth \addtolength{\leftmargin}{\parindent}
 \rightmargin=\parindent
 \raggedright
}
\item[\texttt{subroutine initopt (}]
\texttt{pname,} \texttt{n,} \texttt{mi,} \texttt{me,} \texttt{x,}
\texttt{lx,} \texttt{ux,} \texttt{dxmin,} \texttt{li,} \texttt{ui,}
\texttt{dcimin,} \texttt{inf,} \texttt{tolopt,} \texttt{simcap,}
\texttt{info,} \texttt{izs,} \texttt{rzs,} \texttt{dzs)}
\end{list}
\end{minipage}
$}
\end{quote}
\begin{list}{}{\topsep=0pt\partopsep=0pt\parsep=0.5\parskip\itemsep=0.5\parskip
 \listparindent=\parindent
 \labelwidth=\parindent
 \labelsep=1ex
 \addtolength{\leftmargin}{-1ex}
}
\item[\texttt{pname} (O):]
character string of length 132, giving the name of the problem.

\item[\texttt{n} (I), \texttt{mi} (I), \texttt{me} (I):]
dimensions of the problem. Their meaning is given in the description of
\texttt{dimopt}.

\item[\texttt{x} (O):]
\texttt{double precision} array of dimension $n$, providing a starting
point for the optimization solver.

\item[\texttt{lx} (O), \texttt{ux} (O):]
\texttt{double precision} array of dimension $n$, providing the bounds
on the variable $x$. In other words, $x_i$ is required to satisfy
$\texttt{lx(}i\texttt{)}\leq x_i\leq \texttt{ux(}i\texttt{)}$, for
$i=1,\ldots,n$. The lower (resp.\ upper) bound
\texttt{lx(}$i$\texttt{)} (resp.\ \texttt{ux(}$i$\texttt{)}) is set to
\texttt{-inf} (resp.\ \texttt{inf}) is the bound does not exist; see
below for the meaning of \texttt{inf}. Therefore, the arrays
\texttt{lx} and \texttt{ux} must have been declared in the calling
program with dimension~$n$, even if a solver not dealing with bound
constraints is intended to be used.

\item[\texttt{dxmin} (O):]
\texttt{double precision} variable, providing the resolution in~$x$
for the $l_\infty$ norm: two points whose distance in $\RRn$ for the
{\em sup\/}-norm is less than \texttt{dxmin} can be considered as
indistinguishable. This data can be used in line-search or
trust-region. It is also useful to detect bounds that are active up to
that precision.

\item[\texttt{li} (O), \texttt{ui} (O):]
\texttt{double precision} array of dimension $\texttt{mi}:=m_I$,
providing the bounds on the constraint values $c_I(x)$. In other words,
$c_i(x)$ is required to satisfy $\texttt{li(}i\texttt{)}\leq c_i(x)\leq
\texttt{ui(}i\texttt{)}$, for $i=1,\ldots,m_I$.

\item[\texttt{dcimin} (O):]
\texttt{double precision} variable, providing the resolution in~$c_I$
for the $l_\infty$ norm: two inequality constraint values whose
distance in $\RRmI$ for the {\em sup\/}-norm is less than
\texttt{dcimin} can be considered as indistinguishable. This data can
be useful to detect inequality constraints that are active up to that
precision.

\item[\texttt{inf} (O):]
\texttt{double precision} variable, specifying what is the infinite
value for the bounds on $x$ and $c_I(x)$. In other words, when
$\texttt{lx(}i\texttt{)}\leq-\texttt{inf}$ (resp.\
$\texttt{li(}i\texttt{)}\leq-\texttt{inf}$), there is no lower bound on
$x_i$ (resp.\ $c_i(x)$). A similar convention is adopted for the upper
bounds.

\item[\texttt{tolopt} (O):]
\texttt{double precision} array of dimension $4$, providing the
tolerances on optimality that a pair $(x,\lm)$ must satisfied in order
to be considered as a solution to the problem. More specifically, the
pair $(x,\lm)$ can be considered as a satisfiable KKT point if
\begin{eqnarray*}
\|\nabla_x\ell(x,\lm)\|_\infty&\leq&\texttt{tolopt(1)}\\
\|c(x)^\#\|_\infty&\leq&\texttt{tolopt(2)}\\
\|\sgn_x(\lm)\|_\infty&\leq&\texttt{tolopt(3)},
\end{eqnarray*}
where $\sgn_x(\lm)\in\RRm$ is defined as follows
$$
(\sgn_x(\lm))_i=
\begin{accolade}{ll}
\lm_i^+ & \mbox{if $i\in B\union I$ and
$x_i\notin[l_i+\texttt{tolopt(2)},+\infty[$}\\
\lm_i & \mbox{if $i\in B\union I$ and
$x_i\in[l_i+\texttt{tolopt(2)},u_i-\texttt{tolopt(2)}]$}\\
\lm_i^- & \mbox{if $i\in B\union I$ and $x_i\notin
\hspace{10.67ex}
]-\infty,u_i-\texttt{tolopt(2)}]$}\\
0 & \mbox{if $i\in E$}.
\end{accolade}
$$
The tolerance \texttt{tolopt(4)} is aimed at being used by minimization
software for {\em nonsmooth functions\/}\index{function!nonsmooth} and
provides a tolerance of the duality gap. This way of checking
optimality will probably be improved in a future version of the
collection, in the light of~\cite{dolan-more-munson-2006}.

\item[\texttt{simcap} (O):]
\texttt{integer} array of dimension 4. It specifies the simulator
capabilities. A negative values means that the related function is not
present or that the capability is not considered by the simulator.
\begin{list}{}{\topsep=0.5ex\partopsep=0.5ex\parsep=0.0ex\itemsep=0.0ex
 \listparindent=\parindent
 \settowidth{\labelwidth}{$\texttt{simcap(1)}=9$:}
 \labelsep=1ex
 \leftmargin=\labelwidth
 \addtolength{\leftmargin}{1ex}
}
\item[$\texttt{simcap(1)}<0$]
the simulator cannot evaluate the cost-function~$f$; it may be assumed
then that this one is constant (or zero), so that the problem is a
feasibility one;
\item[$=0$]
the simulator can evaluate the cost-function~$f$;
\item[$=1$]
the cost-function~$f$ is {\em nonsmooth\/}\index{function!nonsmooth}
(this is the only place where this property of the problem can be
detected) and the simulator can evaluate~$f$ and a subgradient~$g$;
\item[$=2$]
the simulator can evaluate the cost-function~$f$ and its gradient~$g$;

\item[$\texttt{simcap(2)}<0$]
the simulator cannot evaluate the inequality constraint function $c_I$;
this is normally because there is no inequality constraints;
\item[$=0$]
the simulator can evaluate~$c_I$;
\item[$=1$]
the simulator can evaluate~$c_I$ and its Jacobian~$c'_I$;

\item[$\texttt{simcap(3)}<0$]
the simulator cannot evaluate the equality constraint function~$c_E$;
this is normally because there is no equality constraints;
\item[$=0$]
the simulator can evaluate~$c_E$;
\item[$=1$]
the simulator can evaluate~$c_E$ and its Jacobian~$c'_E$;

\item[$\texttt{simcap(4)}<0$]
the simulator cannot evaluate $Hv$, the product of the Hessian of the
Lagrangian $H:=\nabla_{xx}^2\ell(x,\lm)$ times a vector~$v$;
\item[$=1$]
the simulator can evaluate a product $Hv$;
\item[$=2$]
the simulator can evaluate the $H$.
\end{list}
\item[\texttt{info} (O):]
\texttt{integer} variable. If negative ($<0$), \texttt{\SOLV\_}\abr
\texttt{modulopt\_}\abr \texttt{main} should consider that the
initialization of the problem by \texttt{initopt} has failed and should
stop.
\item[\texttt{izs}, \texttt{rzs}, \texttt{dzs} (O):]
\texttt{integer}, \texttt{real}, and \texttt{double precision} arrays
that \texttt{initopt} should initialize. These variables are made
available to the Modulopt problem. Their dimensions have been provided
on return from \texttt{dimopt} and they should have been allocated by
the main program \texttt{\SOLV\_}\abr \texttt{modulopt\_}\abr
\texttt{main} associated with some code \SOLV.
\end{list}

\index{subroutine!initopt@\texttt{initopt}|)}

\subsubsection*{The subroutine \texttt{simulopt}}

\index{subroutine!simulopt@\texttt{simulopt}|(}

The subroutine \texttt{simulopt} is the simulator of the problem. It
can be called by \texttt{\SOLV\_}\abr \texttt{modulopt\_}\abr
\texttt{main}, before calling \SOLV. It is also called by the latter to
have information (function and their derivatives) on the problem to
solve. In \texttt{Fortran 95}, it has the following calling structure:
\begin{quote}
\fbox{$
\setlength{\minipagewidth}{\linewidth}
\addtolength{\minipagewidth}{-6.8pt}
\begin{minipage}[t]{\minipagewidth}
\begin{list}{}{\topsep=0pt\partopsep=0pt\parsep=0.5\parskip\itemsep=0.5\parskip
 \settowidth{\labelwidth}{\texttt{subroutine simulopt (}}
 \labelsep=0ex
 \leftmargin=\labelwidth \addtolength{\leftmargin}{\parindent}
 \rightmargin=\parindent
 \raggedright
}
\item[\texttt{subroutine simulopt (}]
\texttt{indic,} \texttt{n,} \texttt{mi,} \texttt{me,} \texttt{x,}
\texttt{lm,} \texttt{f,} \texttt{ci,} \texttt{ce,} \texttt{g,}
\texttt{ai,} \texttt{ae,} \texttt{v,} \texttt{hlv,} \texttt{hl,}
\texttt{izs,} \texttt{rzs,} \texttt{dzs)}
\end{list}
\end{minipage}
$}
\end{quote}
\begin{list}{}{\topsep=0pt\partopsep=0pt\parsep=0.5\parskip\itemsep=0.5\parskip
 \listparindent=\parindent
 \labelwidth=\parindent
 \labelsep=1ex
 \addtolength{\leftmargin}{-1ex}
}
\item[\texttt{indic} (IO):]
\texttt{integer} variable monitoring the communication between the
solver and the simulator. The simulator \texttt{simulopt} recognizes
the following values of \texttt{indic}.
\begin{list}{}{\topsep=0.5ex\partopsep=0.5ex\parsep=0.0ex\itemsep=0.0ex
 \listparindent=\parindent
 \settowidth{\labelwidth}{$\geq 9$:}
 \labelsep=1ex
 \leftmargin=\labelwidth
 \addtolength{\leftmargin}{1ex}
}
\item[$=1$:]
The simulator can do anything except changing the value of the
arguments of \texttt{simulopt}. Typically it prints some information on
the screen, in a file, or on a plotter. Some solver calls the simulator
with this value of \texttt{indic} at each iteration.

\item[$=2$:]
The simulator is asked to compute the value of the functions
$\texttt{f}=f(x)\in\RR$ (cost function), $\texttt{ci}=c_I(x)\in\RRmI$
(inequality constraints), and $\texttt{ce}=c_E(x)\in\RRmE$ (equality
constraints) at a given point~$x$.

\item[$=3$:]
The simulator is asked to compute $\texttt{g}=\nabla f(x)\in\RRn$
(gradient of $f$ at $x$ for the Euclidean scalar product),
$\texttt{ai}=c_I'(x)$ ($m_I\times n$ Jacobian matrix of $c_I$ at~$x$,
hence the $(i,j)$ entry of \texttt{ai} must be the partial derivative
${\partial c_i}/{\partial x_j}$ evaluated at $x$), and
$\texttt{ae}=c_E'(x)$ ($m_E\times n$ Jacobian matrix of $c_E$ at~$x$).

\item[$=4$:]
The simulator is asked to compute $\texttt{f}=f(x)$,
$\texttt{ci}=c_I(x)$, and $\texttt{ce}=c_E(x)$ at a given point $x$, as
well as the gradient $\texttt{g}=\nabla f(x)\in\RRn$,
$\texttt{ai}=c_I'(x)$, and $\texttt{ae}=c_E'(x)$.

\item[$= 5$:]
The simulator is asked to prepare for subsequent computations of
products $Hv$, where $H:= \nabla^2_{xx}\ell(x,\lm)$ is the Hessian of
the Lagrangian at $(x,\lm)$ (see \eqref{lagrangian}) and $v$ will be an
arbitrary vector (see $\texttt{indic}=6$ below). In some cases, it is
convenient to compute the full Hessian of the Lagrangian $H$ when
$\texttt{indic}=5$. In other cases, computing $H$ is too expensive and
nothing has to be done when \texttt{simul} is called with
$\texttt{indic}=5$.

\item[$= 6$:]
The simulator is asked to compute a product $Hv$, where
$H:=\nabla_{xx}^2\ell(x,\lm)$ is the Hessian of the Lagrangian at the
point $(x,\lm)$ and $v\in\RRn$ is an arbitrary vector. Note that it is
not necessary to evaluate the whole matrix $H$ to compute a product
$Hv$; indeed, $Hv$ is also the directional derivative of the gradient
of the Lagrangian in the direction $v$:
$$
Hv=
\left.\frac{d}{dt}\Bigl(\nabla_x\ell(x+tv,\lm)\Bigr)\right|_{t=0}.
$$

\item[$= 7$:]
The simulator is asked to compute the Hessian of the Lagrangian
$H:=\nabla_{xx}^2\ell(x,\lm)$ at the point $(x,\lm)$.
\end{list}

\noindent
On the other hand, the simulator \texttt{simulopt} can also send a
message to the solver, by giving to \texttt{indic} one of the following
values.
\begin{list}{}{%
 \topsep=0pt\partopsep=0pt\parsep=0.5\parskip\itemsep=0.5\parskip
 \listparindent=\parindent
 \settowidth{\labelwidth}{$=-9$:}
 \labelsep=1ex
 \leftmargin=\labelwidth
 \addtolength{\leftmargin}{1ex}
}
\item[$\geq\hfill 0$:]
normal call; the required computation has been done.

\item[$=-1$:]
by this value, the simulator tells the solver that it is impossible or
undesirable to do the calculation at the point $x$ given by the solver.
The reaction of the solver will vary from one solver to the other.


\item[$=-2$:]
the simulator asks the solver to stop, for example because some events
that the solver cannot understand (not in the field of optimization)
has occurred.

\end{list}

\item[\texttt{n} (I), \texttt{mi} (I), \texttt{me} (I):]
dimensions of the problem. Their meaning is given in the description of
\texttt{dimopt}.
\item[\texttt{x} (I):]
\texttt{double precision} array of dimension $n$, providing the point at
which the simulator has to evaluate functions and derivatives.
\item[\texttt{lm} (I):]
\texttt{double precision} array of dimension $m$, providing the current
value of the dual variable $\lambda$. This one determines, with $x$,
the primal-dual variables at which the simulator has to evaluate the
Hessian of the Lagrangian or the product of this Hessian with a vector
(this depends on the value of \texttt{indic}).
\item[\texttt{f} (O):]
\texttt{double precision} variable, providing the cost function value
$f(x)$ if $\texttt{indic}\abr =2$ or $4$ on entry.
\item[\texttt{ci} (O):]
\texttt{double precision} array of dimension $m_I$, providing the
inequality constraint value $c_I(x)$ if $\texttt{indic}=2$ or $4$ on
entry.

\item[\texttt{ce} (O):]
\texttt{double precision} array of dimension $m_E$, providing the
equality constraint value $c_E(x)$ if $\texttt{indic}=2$ or $4$ on
entry.

\item[\texttt{g} (O):]
\texttt{double precision} array of dimension $n$, providing the
gradient of the cost function $\nabla f(x)$ if $\texttt{indic}=3$ or
$4$ on entry.

\item[\texttt{ai} (O):]
\texttt{double precision} array of dimension $m_I\times n$, providing
the Jacobian matrix of the inequality constraint function $c'_I(x)$ if
$\texttt{indic}=3$ or $4$ on entry.

\item[\texttt{ae} (O):]
\texttt{double precision} array of dimension $m_I\times n$, providing
the Jacobian matrix of the equality constraint function $c'_E(x)$ if
$\texttt{indic}=3$ or $4$ on entry.

\item[\texttt{v} (I):]
\texttt{double precision} array of dimension $n$, providing the vector
$v$ that multiplies the Hessian of the Lagrangian if $\texttt{indic}=6$
on entry.

\item[\texttt{hlv} (O):]
\texttt{double precision} array of dimension $n$, providing the product
$Hv$ of the Hessian of the Lagrangian $H$ with a vector $v$ if
$\texttt{indic}=6$ on entry.

\item[\texttt{hl} (O):]
\texttt{double precision} array of dimension $(n,n)$, providing the
Hessian of the Lagrangian $H$ if $\texttt{indic}=7$ on entry.

\item[\texttt{izs}, \texttt{rzs}, \texttt{dzs} (IO):]
\texttt{integer}, \texttt{real}, and \texttt{double precision} arrays
that \texttt{simulopt} can use and modify. These variables are made
available to the Modulopt problem. Their dimensions have been provided
on return from \texttt{dimopt} and they should have been allocated by
the main program \texttt{\SOLV\_}\abr \texttt{modulopt\_}\abr
\texttt{main} associated with some code \SOLV.
\end{list}

\index{subroutine!simulopt@\texttt{simulopt}|)}

\subsubsection*{The subroutine \texttt{postopt}}

\index{subroutine!postopt@\texttt{postopt}|(}

The subroutine \texttt{postopt} is normally called by the main program
\texttt{\SOLV\_}\abr \texttt{modulopt\_}\abr \texttt{main} to allow the
problem to provide a post-optimal analysis. Some problems will take
advantage of this opportunity, but most of them won't (they will
provide a subroutine with en empty body). The most trivial operation
that can be done in this subroutine is to print the solution on the
screen. Another possibility is to check second order optimality. The
flexibility offered by this subroutine will allow the user of \libopt\
to make other job than comparing the effect of using various solvers on
his/her problem.

In \texttt{Fortran 95}, \texttt{postopt} has the following calling
structure:
\begin{quote}
\fbox{$
\setlength{\minipagewidth}{\linewidth}
\addtolength{\minipagewidth}{-6.8pt}
\begin{minipage}[t]{\minipagewidth}
\begin{list}{}{\topsep=0pt\partopsep=0pt\parsep=0.5\parskip\itemsep=0.5\parskip
 \settowidth{\labelwidth}{\texttt{subroutine postopt (}}
 \labelsep=0ex
 \leftmargin=\labelwidth \addtolength{\leftmargin}{\parindent}
 \rightmargin=\parindent
 \raggedright
}
\item[\texttt{subroutine postopt (}]
\texttt{n,} \texttt{mi,} \texttt{me,} \texttt{x,} \texttt{lm,}
\texttt{f,} \texttt{ci,} \texttt{ce,} \texttt{g,} \texttt{ai,}
\texttt{ae,} \texttt{hl,} \texttt{izs,} \texttt{rzs,} \texttt{dzs)}
\end{list}
\end{minipage}
$}
\end{quote}
\begin{list}{}{\topsep=0pt\partopsep=0pt\parsep=0.5\parskip\itemsep=0.5\parskip
 \listparindent=\parindent
 \labelwidth=\parindent
 \labelsep=1ex
 \addtolength{\leftmargin}{-1ex}
}
\item[\texttt{n} (I), \texttt{mi} (I), \texttt{me} (I):]
dimensions of the problem. Their meaning is given in the description of
\texttt{dimopt}.

\item[\texttt{x} (IO), \texttt{lm} (IO):]
\texttt{double precision} arrays of dimension $n$ and $m$ respectively.
They provide the primal ($\texttt{x}=x$) and dual
($\texttt{lm}=\lambda$) variables determined by the solver. They may be
modified, since \libopt\ will no longer use them.

\item[\texttt{f} (IO)\hspace{-1ex}], \texttt{ci} (IO), \texttt{ce}
(IO), \texttt{g} (IO), \texttt{ai} (IO), \texttt{ae} (IO), \texttt{hl}
(IO): variables providing the value of $f(x)$, $c_I(x)$, $c_E(x)$,
$g(x)$, $A_I(x)$, $A_E(x)$, and $\nabla_{xx}^2\ell(x,\lm)$ found by the
last call to \texttt{simulopt} (hence the actual values depend on the
capabilities of the simulator and the design of the solver). See the
description of \texttt{simulopt} for the type and dimension of these
variables. These may be modified, since \libopt\ will no longer use
them.

\item[\texttt{izs}, \texttt{rzs}, \texttt{dzs} (IO):]
\texttt{integer}, \texttt{real}, and \texttt{double precision} arrays
that \texttt{postopt} can use and modify. These variables are made
available to the Modulopt problem. Their dimensions have been provided
on return from \texttt{dimopt} and they should have been allocated by
the main program \texttt{\SOLV\_}\abr \texttt{modulopt\_}\abr
\texttt{main} associated with some code \SOLV.
\end{list}

\index{subroutine!postopt@\texttt{postopt}|)}

\subsubsection*{The subroutines \texttt{inprodopt}, \texttt{ctonbopt},
and \texttt{ctcabopt}}

\index{subroutine!inprodopt@\texttt{inprodopt}|(}

Some unconstrained optimization solvers can deal with gradients that
are associated with an inner product, say $(x,y)\mapsto\<x,y>$,
different from the Euclidean inner product $(x,y)\mapsto x\T y$. Such
an inner product is a way of rescaling the problem. These solvers must
be informed of this inner product and this is the role of the
subroutine \texttt{inprodopt}. We describe the structure of the
subroutine in \texttt{Fortran 95}.
\begin{quote}
\fbox{$
\setlength{\minipagewidth}{\linewidth}
\addtolength{\minipagewidth}{-6.8pt}
\begin{minipage}[t]{\minipagewidth}
\begin{list}{}{\topsep=0pt\partopsep=0pt\parsep=0.5\parskip\itemsep=0.5\parskip
 \settowidth{\labelwidth}{\texttt{subroutine inprodopt (}}
 \labelsep=0ex
 \leftmargin=\labelwidth \addtolength{\leftmargin}{\parindent}
 \rightmargin=\parindent
 \raggedright
}
\item[\texttt{subroutine inprodopt (}]
\texttt{n,} \texttt{v1,} \texttt{v2,} \texttt{ip,} \texttt{izs,}
\texttt{rzs,} \texttt{dzs)}
\end{list}
\end{minipage}
$}
\end{quote}
\begin{list}{}{\topsep=0pt\partopsep=0pt\parsep=0.5\parskip\itemsep=0.5\parskip
 \listparindent=\parindent
 \labelwidth=\parindent
 \labelsep=1ex
 \addtolength{\leftmargin}{-1ex}
}
\item[\texttt{n} (I):]
dimension of the vectors whose inner product is going to be taken.

\item[\texttt{v1} (I), \texttt{v2} (I):]
\texttt{double precision} arrays of dimension $n$.
These are the vectors whose inner product is desired.

\item[\texttt{ip} (O):]
\texttt{double precision} variable representing the inner product of
\texttt{v1} and \texttt{v2}.

\item[\texttt{izs}, \texttt{rzs}, \texttt{dzs} (IO):]
\texttt{integer}, \texttt{real}, and \texttt{double precision} arrays
that \texttt{postopt} can use and modify. These variables are made
available to the Modulopt problem. Their dimensions have been provided
on return from \texttt{dimopt} and they should have been allocated by
the main program \texttt{\SOLV\_}\abr \texttt{modulopt\_}\abr
\texttt{main} associated with some code \SOLV.
\end{list}

\index{subroutine!inprodopt@\texttt{inprodopt}|)}

\medskip
Some unconstrained optimization solvers not only need the inner product
subroutine \texttt{inprodopt} but also subroutines that make a change
of coordinates from the {\em canonical orthogonal basis\/} of $\RRn$ to
some {\em orthogonal basis\/}\index{basis!orthonormal} for the inner
product $\<\cdot,\cdot>$. The {\em canonical orthogonal
basis\/}\index{basis!canonical} of $\RRn$ is the set of vectors
$\{\hat{e}^i\}_{i=1}^n$, where the $j$th component of $\hat{e}^i$ is
equal to $\delta_{ij}$ (the {\em Kronecker symbol}\index{Kronecker
symbol}, which is one when $i=j$ and zero otherwise). If a vector is
written $\sum_ix_i\hat{e}^i$ in the canonical basis and $\sum_iy_ie^i$
in the considered orthogonal basis, the subroutine \texttt{ctonbopt}
gives the coordinates $y:=(y_1,\ldots,y_n)$ from $x:=(x_1,\ldots,x_n)$
and the subroutine \texttt{ctcabopt} gives the coordinates $x$
from~$y$.

For example, suppose that
$$
\<u,v> = u\T M\T Mv,
$$
where $M$ is a nonsingular $n\times n$ matrix, such that a linear
system with the matrix $M$ is easy to solve (for example $M$ could be
triangular). One can take $e^i = M^{-1}\hat{e}^i$, for $1 \leq i \leq
n$, since then $\<e^i,e^j>=(e^i)\T M\T
Me^j=(\hat{e}^i)\T\hat{e}^j=\delta_{ij}$. Knowing the coordinates
$x:=(x_1,\ldots,x_n)$ of a vector in the canonical basis, its
coordinates $y:=(y_1,\ldots,y_n)$ in the basis $\{e^i\}_{i=1}^n$ can be
computed by
$$
y_j
=\<\sum_ix_i\hat{e}^i,e^j>
=\sum_ix_i\<\hat{e}^i,M^{-1}\hat{e}^j>
=\sum_ix_i(\hat{e}^i)\T M\T \hat{e}^j
=(Mx)_j.
$$
We have shown that $y=Mx$. In that example, the subroutine
\texttt{ctonbopt} will compute $y=Mx$ knowing $x$, while the subroutine
\texttt{ctcabopt} will compute $x=M^{-1}y$ knowing~$y$.

Here is the description of the subroutines \texttt{ctonbopt} and
\texttt{ctcabopt} in \texttt{Fortran}~\texttt{95}. The variables
$\texttt{x}=x$ and $\texttt{y}=y$ have the same meaning as in the
discussion above. The parameters \texttt{izs}, \texttt{rzs}, and
\texttt{dzs} have the same meaning as in the subropuitine
\texttt{inprodopt}.
\begin{quote}
\index{subroutine!ctonbopt@\texttt{ctonbopt}|(}
\fbox{$
\setlength{\minipagewidth}{\linewidth}
\addtolength{\minipagewidth}{-6.8pt}
\begin{minipage}[t]{\minipagewidth}
\begin{list}{}{\topsep=0pt\partopsep=0pt\parsep=0.5\parskip\itemsep=0.5\parskip
 \settowidth{\labelwidth}{\texttt{subroutine ctonbopt (}}
 \labelsep=0ex
 \leftmargin=\labelwidth \addtolength{\leftmargin}{\parindent}
 \rightmargin=\parindent
 \raggedright
}
\item[\texttt{subroutine ctonbopt (}]
\texttt{n,} \texttt{x,} \texttt{y,} \texttt{izs,} \texttt{rzs,}
\texttt{dzs)}
\end{list}
\end{minipage}
$}
\index{subroutine!ctonbopt@\texttt{ctonbopt}|)}
\\[1ex]
\index{subroutine!ctcabopt@\texttt{ctcabopt}|(}
\fbox{$
\setlength{\minipagewidth}{\linewidth}
\addtolength{\minipagewidth}{-6.8pt}
\begin{minipage}[t]{\minipagewidth}
\begin{list}{}{\topsep=0pt\partopsep=0pt\parsep=0.5\parskip\itemsep=0.5\parskip
 \settowidth{\labelwidth}{\texttt{subroutine ctcabopt (}}
 \labelsep=0ex
 \leftmargin=\labelwidth \addtolength{\leftmargin}{\parindent}
 \rightmargin=\parindent
 \raggedright
}
\item[\texttt{subroutine ctcabopt (}]
\texttt{n,} \texttt{y,} \texttt{x,} \texttt{izs,} \texttt{rzs,}
\texttt{dzs)}
\end{list}
\end{minipage}
$}
\index{subroutine!ctcabopt@\texttt{ctcabopt}|)}
\end{quote}
Of course, if the inner product implemented in \texttt{inprodopt} is
the Euclidean inner product, \texttt{ctonbopt} will just copy
\texttt{y} into \texttt{x}, while \texttt{ctcabopt} will just copy
\texttt{x} into \texttt{y}.

\subsection{The makefile describing how to run the Modulopt problems}
\label{s:makefile-test-problem}

In the framework introduced in section~\ref{s:new-problem}, the target
\PROB\ of the makefile
\begin{quote}
\texttt{%
\$LIBOPT\_DIR/collections/modulopt/probs/\textit{prob}/Makefile}
\end{quote}
has to specify how to get in the working directory the archive
\texttt{\textit{prob}.a} containing the code of
\texttt{dimopt}, \texttt{initopt}, \texttt{simulopt}, \texttt{postopt},
\texttt{inprodopt}, \texttt{ctonbopt}, and \texttt{ctcabopt}. It must
also create in the working directory symbolic links to the data files
(if any) useful for the execution of the problem. Note that this
makefile is launched from the working directory, so that `\texttt{.}'
in the makefile refers to that directory.

The target \texttt{\textit{prob}\_}\abr \texttt{clean} of
the same makefile, should also clean up the working directory from the
files related to the problem just solved. This makefile is also
launched from the working directory.

One must keep in mind that the makefile must be organized in such a way
that no file is generated in the problem directory
\texttt{\$LIBOPT\_DIR/}\abr \texttt{collections/}\abr
\texttt{modulopt/}\abr \texttt{probs/}\abr \PROB. Only the working
directory can be modified, so that several users can use the Libopt
environment at the same time.

\section{Making a solver able to solve Modulopt prob\-lems}
\label{s:preparing-solver}

In this section, we consider the case when it is desirable to make a
solver of optimization problems, installed in the Libopt environment,
able to solve problems from the Modulopt collection. Let
\begin{quote}
\SOLV\index{solv@\SOLV}
\end{quote}
be the name of the considered solver. If the solver does not already
exist in the Libopt environment, the name \SOLV\ has to be
added to the file
\begin{quote}
\texttt{%
\$LIBOPT\_DIR/solvers/solvers.lst},
\end{quote}
which contains the list of solvers of the Libopt environment,
and the following directories must be created
\begin{quote}
\texttt{%
\$LIBOPT\_DIR/solvers/\textit{solv}}\\
\texttt{%
\$LIBOPT\_DIR/solvers/\textit{solv}/bin},
\end{quote}
as well as an empty file
\begin{quote}
\texttt{%
\$LIBOPT\_DIR/solvers/\textit{solv}/collections.lst},
\end{quote}
which contains the list of collections that the solver \SOLV\
can consider (none if this is the first time \SOLV\ is
introduced in the Libopt environment).

We now give the next steps to follow, together with some explanations.

\begin{list}{}{\topsep=0.5ex\partopsep=0.5ex\parsep=0.5ex\itemsep=0.5ex
 \settowidth{\labelwidth}{9.}
 \labelsep=0.5em
 \leftmargin=\labelwidth
 \addtolength{\leftmargin}{\labelsep}
}
\item[1.]
Create the directory
\begin{quote}
\texttt{%
\$LIBOPT\_DIR/solvers/\textit{solv}/modulopt},
\end{quote}
which will contain the programs/scripts to run \SOLV\ on the
Modulopt problems, and add the name \texttt{modulopt} to the list
\begin{quote}
\texttt{%
\$LIBOPT\_DIR/solvers/\textit{solv}/collections.lst},
\end{quote}
which indicates that \SOLV\ can deal with the Modulopt
collection.

\item[2.]
Create the files
\begin{quote}
\texttt{%
\$LIBOPT\_DIR/solvers/\textit{solv}/modulopt/all.lst}\\
\texttt{%
\$LIBOPT\_DIR/solvers/\textit{solv}/modulopt/default.lst}.
\end{quote}
\begin{list}{{\small$\bullet$}}{\topsep=1.0ex\parsep=0.5ex\itemsep=0.5ex
 \settowidth{\labelwidth}{{\small$\bullet$}}
 \labelsep=0.5em
 \leftmargin=\labelwidth
 \addtolength{\leftmargin}{\labelsep}
 \addtolength{\leftmargin}{\parindent}
 \rightmargin=\parindent
}
\item
The first file (\texttt{all.lst}) must list the problems from the
Modulopt collection that \SOLV\ is able to solve or, more
precisely, those for which it has been conceived. It can contain {\em
comments}\index{comment (in \texttt{*.lst} files)}, which start with the
`\texttt{\#}' character and go up to the end of the line. The easiest
way of doing this is to start with a copy of the file
\begin{quote}
\texttt{%
\$LIBOPT\_DIR/collections/modulopt/all.lst},
\end{quote}
which lists all the Modulopt problems, and to remove from the copied
file those problems that do not have the structure expected by
\SOLV. For example, if \SOLV\ is a solver of
unconstrained optimization problems, remove from the copied file
\texttt{all.lst}, all the problems with constraints.\ffootnote{This is
not well done, since it is not easy to know the feature of the problems
just by looking at the list \texttt{all.lst}.} The features of the
Modulopt problems are given in the files
\begin{quote}
\texttt{%
\$LIBOPT\_DIR/collections/modulopt/all.lst},\\
\texttt{%
\$LIBOPT\_DIR/collections/modulopt/probs/PROBLEMS}.
\end{quote}
Note that other lists exist in the directory
\texttt{\$LIBOPT\_DIR/}\abr \texttt{collections/}\abr
\texttt{modulopt}, which might be more appropriate to start with than
the list \texttt{all.lst}.

\item
The second file above (\texttt{default.lst}) can contain any subset of
the problems listed in the first file (\texttt{all.lst}). This file is
used as the default subcollection when no list is specified in the
\texttt{runopt} script. Therefore, it is often a symbolic link to the
first file \texttt{all.lst}, obtained using the Unix/Linux command
\begin{quote}
\texttt{%
ln -s all.lst default.lst}
\end{quote}
in the directory \texttt{\$LIBOPT\_DIR/}\abr \texttt{collections/}\abr
\texttt{modulopt}.
\end{list}

\item[3.]
Create the executable file
\begin{quote}
\texttt{%
\$LIBOPT\_DIR/solvers/\textit{solv}/modulopt/\textit{solv}\_modulopt}.
\end{quote}
This is the script launched by \texttt{runopt} to run \SOLV\ on a
single Modulopt problem. More precisely, this script takes care of the
tasks that can be described using Unix/Linux commands, such as making
symbolic links, executing makefiles, removing files, etc. This is not
an easy program to write from scratch (see~\cite{gilbert-jonsson-2007})
but adapting a script used by another solver is rather straightforward.
For example, one can copy, rename, and modify the Perl script
\begin{quote}
\texttt{%
\$LIBOPT\_DIR/solvers/sqppro/modulopt/sqppro\_modulopt}.
\end{quote}
The only changes to make in the copied and renamed script (hence now
called \texttt{\textit{solv}\_modulopt}) consists in substituting the
two occurrences of the solver name \texttt{sqppro} by \SOLV\ (one
in a comment, another in the definition of the variable
\texttt{\$solv}\-\texttt{name}). That's all!

\item[4.]
Create the main program
\begin{quote}
\texttt{%
\$LIBOPT\_DIR/solvers/\textit{solv}/modulopt/\textit{solv}\_modulopt\_main.f90}.
\end{quote}
This program is very solver dependent and is, with the next step to
which it is linked, the most difficult task to realize. It is the main
program that will be linked with the subroutines describing the problem
from the Modulopt collection selected by the \texttt{runopt} script,
those in the archive \texttt{\textit{prob}.a} (if the selected problem
is \PROB) created by the makefile \texttt{\$LIBOPT\_DIR/}\abr
\texttt{collections/}\abr \texttt{modulopt/}\abr \texttt{probs/}\abr
\PROB\texttt{/}\abr \texttt{Makefile} (see
section~\ref{s:makefile-test-problem}).
The language used to write this main program is arbitrary, provided it
(or its object form generated by some compilor) can be linked with the
object files issued from the compilation of the Fortran files
describing the Modulopt problems (\texttt{dimopt}, \texttt{initopt},
\texttt{simulopt}, \texttt{postopt}, \texttt{inprodopt},
\texttt{ctonbopt}, and \texttt{ctcabopt}).

If Fortran 90/95 is the adopted language, the easiest way to proceed is
to copy and rename the file
\begin{quote}
\texttt{%
\$LIBOPT\_DIR/solvers/sqppro/modulopt/sqppro\_modulopt\_main.f90}
\end{quote}
into the file
\begin{quote}
\texttt{%
\$LIBOPT\_DIR/solvers/\SOLV/modulopt/\SOLV\_modulopt\_main.f90}.
\end{quote}
Since this main program is very solver dependent, its part dealing with
the solver will have to be thoroughly modified. Let us describe the
structure of the program.
\begin{list}{}{\topsep=0.5ex\partopsep=0.5ex\parsep=0.5ex\itemsep=0.5ex
 \settowidth{\labelwidth}{9.9.}
 \labelsep=0.5em
 \leftmargin=\labelwidth
 \addtolength{\leftmargin}{\labelsep}
}
\item[4.1.]
\label{point-4.1}
After the declaration of variables, the program calls the subroutine
\texttt{dimopt}\index{subroutine!dimopt@\texttt{dimopt}}
to get the dimensions of the Modulopt problem that will be selected by
the \texttt{runopt} script. These dimensions are then used to allocate
dimension dependent variables, including \texttt{izs}, \texttt{rzs},
and \texttt{dzs}.

\item[4.2.]
The problem data are then obtained by calling the subroutine
\texttt{initopt}\index{subroutine!initopt@\texttt{initopt}}. This is
the good spot to verify that the features of the problem are compatible
with the solver capabilities, using the variable \texttt{simcap}.

\item[4.3.]
Some optimization solver requires that the simulator be called before
launching the optimization. In this case, this is the good spot for
doing so, by calling
\texttt{simulopt}\index{subroutine!simulopt@\texttt{simulopt}}.

\item[4.4.]
Next, the program calls the optimization solver \SOLV, after
having initialized its arguments and opened relevant files.

\item[4.5.]
Once the optimization has been completed, it is important to write the
\texttt{libopt} line\index{libopt!line}, which summarizes the
performance of the solver \SOLV\ on the currently solved Modulopt
problem. See~\cite{gilbert-jonsson-2007} or the \texttt{libopt} man
page.

\item[4.6.]
It is nice to let the problem do its post-optimal analysis (if any) by
finally calling
\texttt{postopt}\index{subroutine!postopt@\texttt{postopt}}.

\end{list}
Note that a particular solver usually requires a simulator with another
structure than the one of \texttt{simulopt}. Therefore an interface
between \texttt{simulopt} and the simulator required by \SOLV\
should be written and placed in the file \texttt{\textit{solv}\_}\abr
\texttt{modulopt\_}\abr \texttt{main.f90}.

\item[5.]
Create the makefile
\begin{quote}
\texttt{%
\$LIBOPT\_DIR/solvers/\textit{solv}/modulopt/Makefile}.
\end{quote}
The aim of this makefile is to tell the Libopt environment how to link
the solver binary with the object files describing the Modulopt problem
selected by the \texttt{runopt}\index{runopt@\texttt{runopt} (script)}
script. If the latter is \PROB, the corresponding
object files will be at link time in the working directory in the
archive \texttt{\textit{prob}.a} (see section~\ref{s:test-problem}). The
easiest way of doing this is to start with an existing makefile, like
\begin{quote}
\texttt{%
\$LIBOPT\_DIR/solvers/sqppro/modulopt/Makefile}.
\end{quote}
This one will be copied and renamed into the file
\begin{quote}
\texttt{%
\$LIBOPT\_DIR/solvers/\SOLV/modulopt/Makefile}
\end{quote}
and then modified.
\end{list}

You can now try the command
\begin{quote}
\texttt{%
echo "\textit{solv} modulopt \textit{prob}" | runopt -v}
\index{runopt@\texttt{runopt} (script)}
\end{quote}
where the flag \texttt{-v} (verbose) is used to get detailed comments
from the Libopt scripts, which then tell what they actually
do. The flag \texttt{-t} (test mode) can be used instead, if you want
to see what the scripts would do without asking them to do it.

\section{Directories and files}
\label{s:directories-files}

In this section, we list some important directories\index{directory}
and files encountered in this note. Recall that
\texttt{\$LIBOPT\_DIR}\index{environment
variable!\texttt{\$LIBOPT\_DIR}}\index{libopt!environment variable} is
the environment variable that specifies the root
directory\index{directory!libopt root@\texttt{libopt}
root}\index{libopt!root directory} of the Libopt hierarchy. Below,
\SOLV\ is the generic name of a particular solver known to the Libopt
environment.
\begin{list}{{\small$\bullet$}}{\topsep=1.0ex\parsep=0.0ex\itemsep=1.0ex
 \settowidth{\labelwidth}{{\small$\bullet$}}
 \labelsep=0.5em
 \leftmargin=\labelwidth
 \addtolength{\leftmargin}{\labelsep}
 \addtolength{\leftmargin}{0.5\parindent}
}
\item
\texttt{\$LIBOPT\_DIR/collections}\index{directory!collection root}:\\
directory of the collections of problems the Libopt environment can
deal with.
\item
\texttt{\$LIBOPT\_DIR/collections/collections.lst}:\\
list of collections known to and installed into Libopt.
\item
\texttt{\$LIBOPT\_DIR/collections/modulopt}:\\
root directory of the Modulopt collection\index{directory!Modulopt
root} in the Libopt environment.
\item
\texttt{\$LIBOPT\_DIR/collections/modulopt/all.lst}:\\
list of all the problems of the Modulopt collection.
\item
\texttt{\$LIBOPT\_DIR/collections/modulopt/probs}:\\
directory containing one sub-directory for each problem of the Modulopt
collection.
\item
\texttt{\$LIBOPT\_DIR/solvers}:\\
root directory of the solvers\index{directory!solver root} the Libopt
environment can deal with.
\item
\texttt{\$LIBOPT\_DIR/solvers/solvers.lst}:\\
list of solvers known to and installed into Libopt.
\item
\texttt{\$LIBOPT\_DIR/solvers/\SOLV/collections.lst}:\\
list of collections the code \SOLV\ has been prepared to deal with.
\item
\texttt{\$LIBOPT\_DIR/solvers/\SOLV/modulopt}:\\
directory containing the scripts and programs specifying how to run the
code \SOLV\ on problems from the Modulopt collection.
\item
\texttt{\$LIBOPT\_DIR/solvers/\SOLV/modulopt/all.lst}:\\
list of problems from the Modulopt collection, for which the solver
\SOLV\ is designed.
\item
\texttt{\$LIBOPT\_DIR/solvers/\SOLV/modulopt/\SOLV\_modulopt}:\\
Perl script specifying the Unix/Linux commands useful to run the solver
\SOLV\ on a single problem of the Modulopt collection.
\item
\texttt{\$LIBOPT\_DIR/solvers/\SOLV/modulopt/\SOLV\_modulopt\_main.f90}:\\
Fortran 90/95 main program that is used to run \SOLV\ on a Modulopt
problem selected by a \texttt{runopt} script, say \PROB. This
program is linked with the object files (gathered in the archive
\texttt{\PROB.a} in the working directory) describing \PROB.
\end{list}

{\small

\newcommand{\HOME}{$HOME}

\bibliography{\HOME/bibliographies/optimisation/OPTIM,%
  \HOME/bibliographies/algebre/ALGEBRE,%
  \HOME/bibliographies/applications/APPLICATIONS,%
  \HOME/bibliographies/geometrie/GEOM,%
  \HOME/bibliographies/informatique/INFO,%
  \HOME/bibliographies/methodes_numeriques/METHNUM}

\begin{thebibliography}{1}

\bibitem{bongartz-conn-gould-toint-1995}
I.~Bongartz, A.R.Conn, N.I.M. Gould, Ph.L. Toint (1995).
\newblock {CUTE}: Constrained and unconstrained testing environment.
\newblock {\em ACM Transactions on Mathematical Software}, 21, 123--160.

\bibitem{dolan-more-munson-2006}
E.D. Dolan, J.J. Mor\'e, T.S. Munson (2006).
\newblock Optimality measures for performance profiles.
\newblock {\em SIAM Journal on Optimization}, 16, 891--909.

\bibitem{gilbert-jonsson-2007}
J.Ch. Gilbert, X.~Jonsson (2007).
\newblock {LIBOPT} -- {A}n environment for testing solvers on heterogeneous
  collections of problems.
\newblock Technical report, INRIA, BP 105, 78153 Le Chesnay, France.
\newblock (to appear).

\bibitem{gould-orban-toint-2003}
N.~Gould, D.~Orban, Ph.L. Toint (2003).
\newblock {CUTE}r (and {S}if{D}ec), a {C}onstrained and {U}nconstrained
  {T}esting {E}nvironment, revisited.
\newblock {\em ACM Transactions on Mathematical Software}, 29, 373--394.

\bibitem{lemarechal-1980}
C.~Lemar\'{e}chal (1980).
\newblock Using a {M}odulopt minimization code.
\newblock Unpublished technical Note.

\end{thebibliography}

\bibliographystyle{plain_a}
}

\printindex

\end{document}